[1994/12/01]
\documentclass[printscheme]{aomart}
\usepackage{amsmath}
\usepackage{amssymb,amsfonts,latexsym,cancel}
\usepackage{array}
\usepackage{graphicx}
\usepackage{epstopdf}
\usepackage{booktabs,mathptmx}
\usepackage{upgreek}
\usepackage[version=4]{mhchem}
\usepackage{amsthm,empheq}
\usepackage{mathtools}
\usepackage[thinc]{esdiff}
\usepackage{blkarray}
\usepackage{nccmath}
\usepackage{apacite}
\usepackage{setspace}

\usepackage{mathrsfs}
\usepackage{bm}

\usepackage{chemfig}
\usepackage{chemscheme} 

\usepackage{chemmacros}
\chemsetup{modules={scheme}}
\DeclareChemTranslation{scheme}{spanish}{esquema} 
\DeclareChemTranslation{Scheme}{spanish}{Esquema}

\usepackage{tikz}
\usetikzlibrary{tikzmark,calc}

\usepackage{float}
\floatstyle{plaintop}
\restylefloat{table}
\usepackage{tabularx}

\makeatletter
\newcommand\mathcircled[1]{%
  \mathpalette\@mathcircled{#1}%
}
\newcommand\@mathcircled[2]{%
  \tikz[baseline=(math.base),red] \node[draw,circle,inner sep=0.5pt] (math) {$\m@th#1#2$};%
}
\makeatother
\makeatletter
\newcommand\mathcircledblue[1]{%
  \mathpalette\@mathcircledblue{#1}%
}
\newcommand\@mathcircledblue[2]{%
  \tikz[baseline=(math.base),blue] \node[draw,circle,inner sep=0.5pt] (math) {$\m@th#1#2$};%
}
\makeatother
\makeatletter
\newcommand\mathcirclednaranja[1]{%
  \mathpalette\@mathcirclednaranja{#1}%
}
\newcommand\@mathcirclednaranja[2]{%
  \tikz[baseline=(math.base),orange] \node[draw,circle,inner sep=0.5pt] (math) {$\m@th#1#2$};%
}
\makeatother
\makeatletter
\newcommand\mathcircledcyan[1]{%
  \mathpalette\@mathcircledcyan{#1}%
}
\newcommand\@mathcircledcyan[2]{%
  \tikz[baseline=(math.base),cyan] \node[draw,circle,inner sep=0.5pt] (math) {$\m@th#1#2$};%
}
\makeatother
\makeatletter
\newcommand\mathcircledmarron[1]{%
  \mathpalette\@mathcircledmarron{#1}%
}
\newcommand\@mathcircledmarron[2]{%
  \tikz[baseline=(math.base),brown] \node[draw,circle,inner sep=0.5pt] (math) {$\m@th#1#2$};%
}
\makeatother
\makeatletter
\newcommand\mathcircledamarillo[1]{%
  \mathpalette\@mathcircledamarillo{#1}%
}
\newcommand\@mathcircledamarillo[2]{%
  \tikz[baseline=(math.base),yellow] \node[draw,circle,inner sep=0.5pt] (math) {$\m@th#1#2$};%
}
\makeatother

\newcommand*\circled[1]{\tikz[baseline=(char.base)]{
            \node[shape=circle,draw,inner sep=2pt] (char) {#1};}}
\newcommand*\squared[1]{\tikz[baseline=(char.base)]{
            \node[shape=rectangle,draw,inner sep=2pt] (char) {#1};}} 
            
\usepackage{csquotes}
\usepackage{parskip} 
\usepackage{blindtext}
\usepackage{lscape}
\usepackage{pdflscape}

\usepackage{hyperref}

\usepackage{afterpage}

\makeatletter 
\fancypagestyle{firstpage}{%
  \fancyhf{}%
  \if@aom@manuscript@mode
    \lhead{\begin{picture}(0,0)%
        \put(-21,-25){\usebox{\@aom@linecount}}%
      \end{picture}}
  \fi
  \chead{}%
   \cfoot{\footnotesize\thepage}}%
\doinumber{}
\makeatother

\copyrightnote{}
\chardef\bslash=`\\ 
\newcommand{\ntt}{\normalfont\ttfamily}

\newtheorem[{}\it]{thm}{Theorem}[section]
\newtheorem{cor}[thm]{Corollary}
\newtheorem{lem}[thm]{Lemma}

\theoremstyle{definition}
\newtheorem{defn}{Definition}[section]

\newtheorem*[{}\it]{notation}{Notation}

\newcommand{\eval}[2][\right]{\relax
  \ifx#1\right\relax \left.\fi#2#1\rvert}

\newcommand{\envert}[1]{\left\lvert#1\right\rvert}

\newcommand{\enVert}[1]{\left\lVert#1\right\rVert}

\title{On the fractal nature of the partition function $p(n)$ and the divisor functions $d_i(n)$}
\author{Romulo Leoncio Cruz Simbron}
\address{Laboratorio de Investigación Fisicoquímica (Labinfis)\\
Universidad Nacional de Ingeniería\\
Lima, Perú}

\address{Laboratorio BigData\\
Centro de Tecnologías de Información y Comunicaciones(CTIC)\\
Universidad Nacional de Ingeniería\\
Lima, Perú}

\address{Instituto de Investigación Ambiental y Microbiano (ININAM) \\
Lima, Perú}

\address{Visiting Scholars \\
Blue Marble Space Institute of Science (BMSIS)\\
Washington, USA}
\email{romulo.cruz.s@uni.pe} 
\fulladdress{Laboratorio de Investigación Fisicoquímica\\
Universidad Nacional de Ingeniería\\
Lima, Perú}
\givenname{Romulo}
\surname{Cruz}


\begin{document}

\begin{abstract}
The partitions of the integers can be expressed exactly in an iterative and closed-form expression. This equation is derived from distributing the partitions of a number in a network that locates each partition in a unique and orderly position. From this representation an iterative equation for the function of the number of divisors was derivated. Also, the number of divisors of a integer can be found from a new function called the trace of the number n, trace(n). As a final preliminary result, using the Bressoud-Subbarao theorem, we obtain a new iterative representation of the sum of divisor function. Using this theorem it is possible to derive an iterative equation for any divisor function and all their networks representation will exhibits a self-similarity behavior. We must then conclude that the intricate nature of the divisor functions results from the fractal nature of the partition function described in the present work.
\end{abstract}

\maketitle
\tableofcontents

\section{Partition function p(n)}
The partition p (n) of a natural number n represent the number of ways in which n can be written as an sum of natural numbers without taking into account the order of the additions \cite{ono2000distribution,andrews1998theory,hardy1918asymptotic,rademacher1938partition}. For example, the partitions of 7 are:\\
\begin{center}
7\\
6 + 1\\
5 + 1 + 1\\
4 + 1 + 1 + 1\\
3 + 1 + 1 + 1 + 1\\
2 + 1 + 1 + 1 + 1 + 1\\
1 + 1 + 1 + 1 + 1 + 1 + 1\\
5 + 2\\
4 + 1 + 2\\
3 + 1 + 1 + 2\\
2 + 1 + 1 + 1 + 2\\
3 + 2 + 2\\
2 + 1 + 2 + 2\\
4 + 3\\
3 + 1 + 3
\end{center}
Since there are 15 ways to express 7, then p (7) is 15. This introduction section is to determine a methodology that allows writing all partitions of a number without repetitions and whose position is unique within each representation. In order to accomplish this task, it will be very important that the methodology followed establishes an order in a certain partition in its application, that is, that a partition can only appear in a single way. We will start this task by defining some important concepts that we will call jumps. When we refer to a partition in general we will refer to the partition of the same number n.\\

\begin{defn}
\textbf{The terms of a partition} are called each of the additions that form a partition of a certain number. These terms are written in a specific order.
\end{defn}
Example:\\
8 5 6 is the partition of the number 19 and 8 is the first term, 5 the second term and 6 the third and last term.\\

\begin{defn}
\textbf{A jump of order 1} is called a transformation of a partition that provides another partition and that consists of taking its first term and dividing it into (a-1) and 1. This jump requires that the term a-1 formed is not less than the last term of the partition.
\end{defn}
Example:\\
5 1 3 $\to$ 4 1 1 3\\
8 5 $\to$ 7 1 5\\

\begin{defn}
\textbf{A jump of order r} is called a transformation of a partition that provides another partition and consists of taking its first term a and partitioning it into (a-r) and r. This jump requires that the a-r term formed is not less than the last term of the partition and that the r term is not greater than the terms in its right.
\end{defn}
Example:\\
r = 3\\
10 4 $\to$ 7 3 4\\
r = 5\\
15 6 $\to$ 10 5 6\\
It follows from this definition that if a partition starts with a single term n, and a first jump of order r is performed. All the jumps that could be applied to this new partition could only have an order r or an order less than r.\\

\begin{defn}
\textbf{A non-ascending sequence of ordered jumps} is one in which jumps are successively applied to a partition, the orders of which can remain constant or go down by one or more unit but never go up.
\end{defn}
Let's start for example with partition 7 and apply an initial jump of order 1. Then the non-ascending sequence of ordered jumps of 7 would be: 7 $\to$ 6 1 $\to$ 5 1 1 $\to$ 4 1 1 1 $\to$ 3 1 1 1 1 $\to$ 2 1 1 1 1 1 $\to$ 1 1 1 1 1 1 1, all jumps of order 1. However, if we start with a jump of order two, the non-ascending sequence that we find would be: 7 $\to$ 5 2 $\to$ 3 2 2, all jumps of order 2, and if we also do jumps of order 1  we get 5 2 $\to$ 4 1 2 $\to$ 3 1 1 2 $\to$ 2 1 1 1 2, and 3 2 2 $\to$ 3 1 2 2, all of them jumps of order 1. We can order this non-ascending sequence of jumps as follows:\\
\\
7\\
\\
5 2\\
4 1 2\\
3 1 1 2\\
2 1 1 1 2\\
\\
3 2 2\\
3 1 2 2\\

\begin{lem}
If a partition $a_{1}\; a_{2}\; a_{3}\; a_{4}\; \cdots\;  a_{m}$ is the result of a non-ascending sequence of jumps starting from a given number n, the sequence $a_{2}\; a_{3}\; a_{4}\; \cdots\; a_{m}\; a_{1}$ is a non-descending sequence of terms. \\
\end{lem}

\begin{proof}
As can be seen from the way in which the terms of a sequence are generated, the term $a_{m}$ is the order of the first jump, the previous term $a_{m-1}$ can only be of the order of $a_{m}$ or less, therefore the term $a_{m- 1}$ is equal to $a_{m}$ or less. Similar reasoning leads us to deduce that the sequence $a_{2}\; a_{3}\; a_{4}\; \cdots \; a_{m}$ is a non-descending sequence of terms. As a general rule of jumps, in each jump the first term can never be less than the last term, then the sequence $a_{2}\; a_{3}\; a_{4}\; \cdots\; a_{m}\; a_{1}$ forms a sequence of non-descending terms.
\end{proof}

\begin{defn}
\textbf{A complete and ordered non-ascending sequence of jumps} is one in which all possible ordered jumps are applied to a number n. 
\end{defn}

For example, for the case of 7 we have the following possible complete ordered sequences.\\
7\\
\\
Initial jump of order 1\\
6 1\\
5 1 1\\
4 1 1 1\\
3 1 1 1 1\\
2 1 1 1 1 1\\
1 1 1 1 1 1 1\\
\\
Initial jump of order 2\\
5 2\\
4 1 2\\
3 1 1 2\\
2 1 1 1 2\\
\\
3 2 2\\
2 1 2 2\\
\\
Initial jump of order 3\\
4 3\\
3 1 3\\
\\

\begin{lem}
The elements of an ordered and complete non-ascending sequence of jumps starting from the number n constitute all the partitions of the number n.
\end{lem}

\begin{proof}
To carry out the proof, we will first demonstrate the uniqueness of the representation of a partition, which involves demonstrating that two representations of a partition whose terms have only been changed in order cannot be present in an ordered and complete non-ascending sequence of jumps. We will also prove that there cannot be partitions with the same representation. This allows us to say that each representation is unique and occupies a certain position within this sequence. Finally, we will prove that given any partition of n, we can place it unequivocally within a complete and ordered non-ascending sequence of jumps starting from n.\\

In the case of having two representations of a partition with only the elements exchanged, the second term must be the same in both representations since it must be the smallest term of all the term present. The third term if it exists must also be the same in both cases. For example, if this third term is u in one representation and v in the other. If u is greater than v in the first representation, v would no longer appear in this representation unless it is equal to u. Since v must necessarily appear in this representation, v could only appear if v equals u, which is contradictory. With the same reasoning it follows that all the terms up to the penultimate term are equal and occupy the same order. In the case of the last term, if in one of the representations this is q and in the other it is r, the first term of this first representation would be r, and vice versa, the first term of the second representation would be q. If q and r are different, in any of the representations the first term will be less than the last, in contradiction with the construction rules of partitions, which, q and r must be the same in both representations. In this way it is shown that there is a single representation for each partition found.\\

We will now prove that the same partition cannot appear twice. To test it, we only have to prove that the path that a partition produces is unique. For that we will prove that the predecessor of a partition is unique. The intermediate term more closest to the left term informs us of the jump it has come from, for example, if this partition is of the form `` a 1 1 ... 1 ... b ... bc ... c. .. '' the order of the last jump made is 1. Going back in this jump gives us the predecessor partition `` a+1 1 ... 1 ... b ... bc ... c ... '' which will be the same predecessor in both cases. If k is the number of 1's present, the two partitions considered should have as their predecessor `` (a + k) b ... b c ... c ... ''. It can be deduced that each partition can only occur in one unique way, that is, following a unique path. This shows us that any partition can be unequivocally located within a complete and ordered non-ascending sequence of jumps.\\

Since if we start the sequence from the number n, any partition we get to will be a partition of n. Therefore, the set of elements of the non-ascending sequence of ordered and completed jumps is contained within the elements that make up the partition of the number n. To prove that these two sets are equal, simply prove that that every partition of a number n is contained within a non-ascending sequence of ordered and completed jumps. To prove this, we only have to prove that given a partition we can reach the number n by inverse of the jumps of a certain order. We must start by ordering the unordered partition by placing the terms in non-descending order. After that, the last term of this ordering would be the first element of the representation that we are using in the present work `` a 1 1 ... 1 ... b ... bc ... c ... '' . This representation is unique for any partition of the number n since the terms keep a non-descending order from the second term to the final term and first. As we have seen in the previous paragraph, a partition of the type `` a 1 1 ... 1 ... b ... bc ... c ... '' can have as its last ancestor the number (a + kb + lc + md ...) where n = a + kl + lc + md. Therefore, whatever the unordered partition, we can order it in a non-descending way as mentioned and we can place it unequivocally within the complete and ordered non-ascending sequence of jumps. This character of complete is required since we must take all the possible jumps. For the above mentioned it is shown that the elements of a non-ascending sequence of jumps contains all the possible partitions of a certain number.
\end{proof}

\subsection{Theorem in the partition of a natural number} 
\begin{figure}
\centering
\includegraphics[scale=0.40]{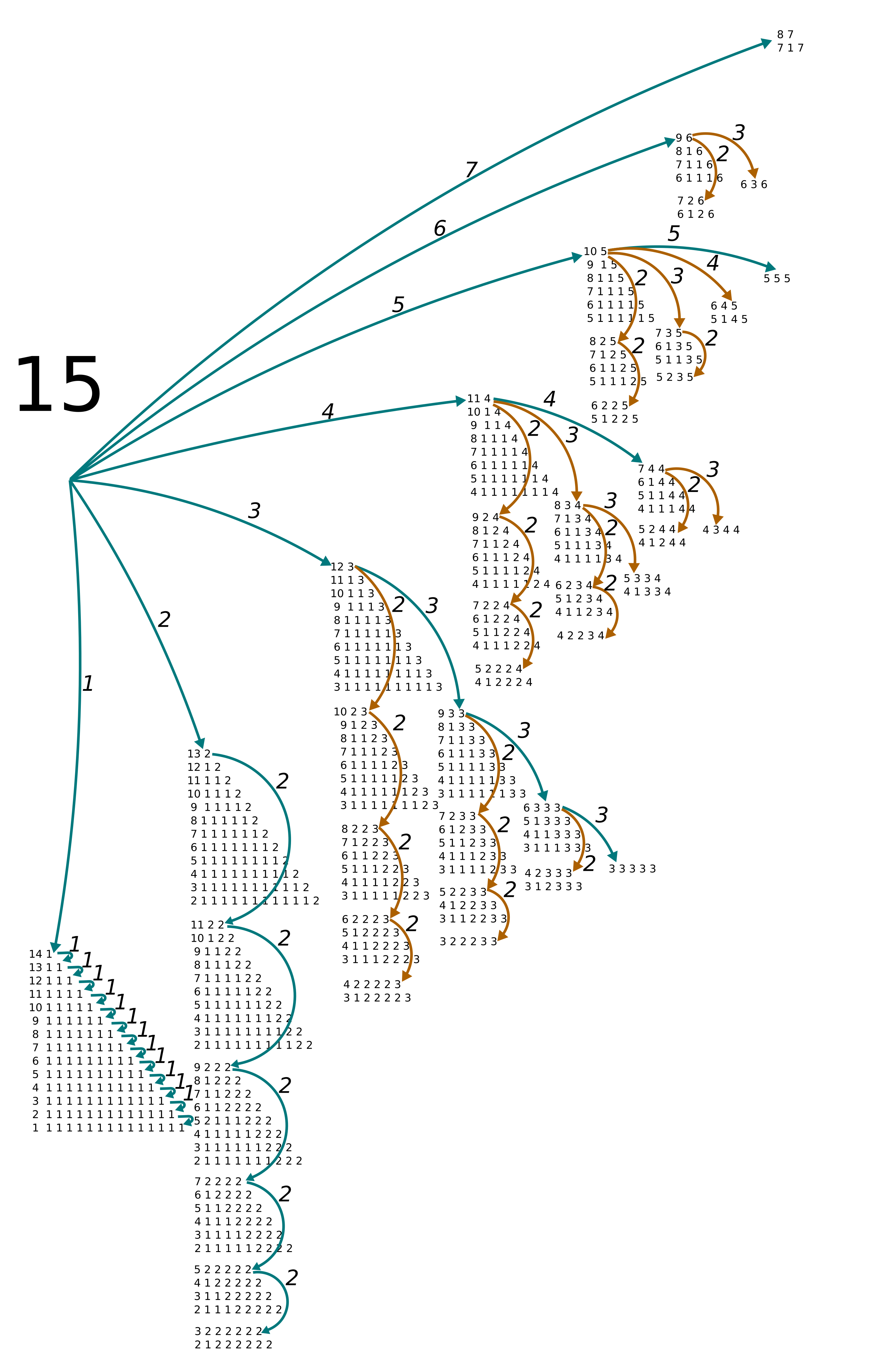}
\caption[]{Layaout of the representation of the integer number 15.}
\label{Secuencia}
\end{figure}

In this section we will detail a form of graphical representation of an ordered and complete non-ascending sequence of jumps. As can be seen from the definition and examples shown in the definition of this type of sequence, it is not a linear sequence but will undoubtedly be a branched representation. In Figure \ref{Secuencia} we show a possible representation of the partitions of the number 15. It can be noted that with this representation all possible non-ascending ordered jumps have been shown. The green lines in the Figure \ref{Secuencia} show the order of the jump performed. These first jumps go from order 1 to order 7. After each jump one can make jumps of an order equal to or less than the jump of the order that preceded it (brown lines). In those groups of jumps that have started with a jump of order greater than 1, the indication of the jumps of order 1 has no longer been placed, but they have been grouped according to the jump of order two from which they come. Following this criterion, the group of partitions that are vertically aligned come from successive jumps of order two. 

\begin{figure}
\centering
\includegraphics[scale=0.80]{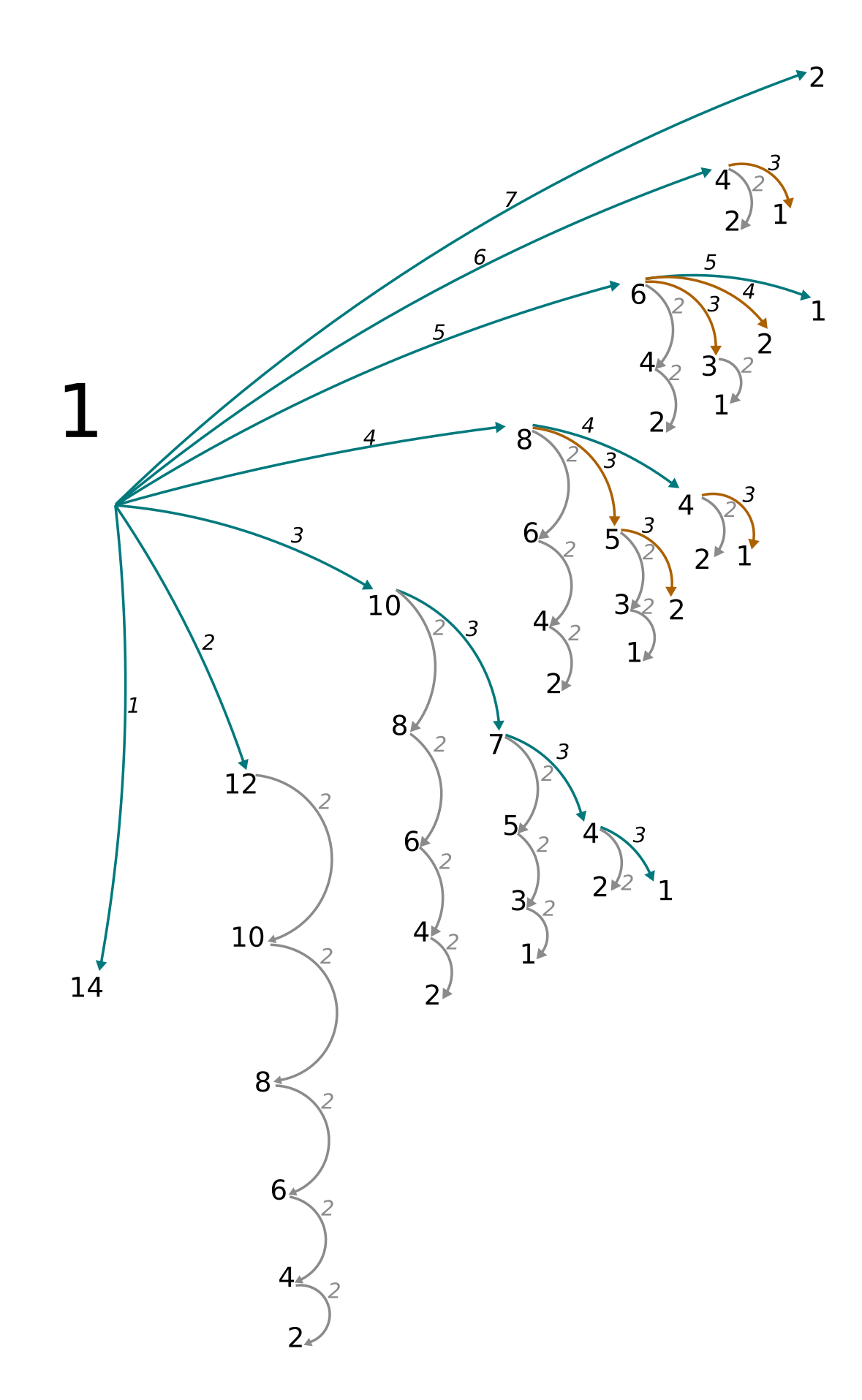}
\caption[]{Layaout of the número of partitions of the number 15.}
\label{Secuencia_numeros}
\end{figure}

In order to elucidate the underlying numerical relationship within this network, what we will do next is to replace the blocks that arise from jumps of order 1 after the initial jumps by the number of partitions they contain. This transformation is shown in Figure \ref{Secuencia_numeros}. We can notice in this figure that descending progressions of ratio 2 occur vertically after each jump. These progressions are complete, that is, they go up to 2 or 1. Could this be a rule? In the following Theorem we show that this is indeed a rule.

\begin{defn}
We will name \textbf{jump set of order 1} to the set of partitions that are the result of the continuous application   of order 1 jumps, including the original partition. 
\end{defn}

For example,  a jump set of order 1 starting from partition 10 4 is:\\
10 4\\
9 1 4\\
8 1 1 4\\
7 1 1 1 4\\
6 1 1 1 1 4\\
5 1 1 1 1 1 4\\
4 1 1 1 1 1 1 4\\
The number of elements in this jump set of order 1 is 7.\\

\begin{lem}
The number of possible jumps of order 1 of any partition that begins with the term a and ends with the term m (`` a b c ... m '' , with $b >1$) is a-m. If we include the original partition, then the number of elements of the jump set of order 1 of the partition `` a b c ... m '' is a-m+1.
\end{lem}

\begin{proof}
In fact, each jump of order 1 decreases the first term by 1 and adds the number 1 as the second term. This type of jump can only be performed until this first term is equal to the first term m, therefore, the number of jumps of order 1 of a partition of type `` a b c ... m '' would be a-m. It follows then that the number of elements in a jump set of order 1 is equal to a-m+1.
\end{proof}

\begin{lem}
If a jump of order $r> 1$ is made from a partition `` a b c ... m '' to a partition `` a-r r b c ... '', the number of elements in the jump set of order 1 of the the partition `` a-r r b c ... m '' decreased by r with respect to the number of elements of the jump set of order 1 of the original partition `` a b c ... m ''.
\end{lem}
\begin{proof}
Since the partition resulting from the jump of order r is `` a-r r b c ... m '', the number of elements of the jump set of order 1 that arise from this partition is a-r-m+1. This number is smaller in r with respect to the number of elements of the jump set of order 1 that arise from the original partition, a-m+1.
\end{proof}

\begin{thm}
If n is a natural number and p(n) represent the number of partitions of the number n, then 

\begin{equation}
\centering
p(n) = \sum_{i=0}^{floor(\frac{n}{2})} \mathscr{P}(i,n+1-2i)
\end{equation}
\begin{equation}
\centering
\mathscr{P}(s,t) =  \sum_{i=1}^{ceil(\frac{t}{s})} \mathscr{P}(s-1,t-s(i-1))
\end{equation}
\begin{equation}
\centering
\mathscr{P}(0,t) = 1
\end{equation}

\end{thm}

\begin{proof}
We can notice that when a first jump of order r is performed, the number of elements of the jump set of order 1 of the resulting partition is n+1-2r. Indeed, if we start from the number n, the partition resulting from a jump of order r is n-r r, the number of elements of the set of jumps of order 1 of this partition is n-r-r + 1, that is, n+1-2r. After these first jumps it is possible to perform all possible non-ascending and ordered jumps and as has been shown in the previous lemmas, after a jump of order s, the number of elements of the jump set of order 1 to which the resulting partition belongs decreases in s with respect to the partition that originated it. Then it follows that the number of partitions of the jump sets 1 resulting from the first jumps follow the peculiar iterative formula.
\begin{equation}
\centering
\mathscr{P}(s,t) =  \sum_{i=1}^{ceil(\frac{t}{s})} \mathscr{P}(s-1,t-s(i-1))
\end{equation}
\begin{equation}
\centering
\mathscr{P}(0,t) = 1
\end{equation}

Where $\mathscr{P}(s,t)$ gives us the number of partitions that can arise after having made an initial jump of order s starting from a partition that originates a set of jumps of order 1 whose elements are in total t. Note that the upper limit of the summation, ceil ($\frac{t}{s}$), comes from the fact that if the jump set of order 1 of the resulting partition is t, then the following partitions that result from successive jumps of order s will be t, t-s , t-2s, t-3s ..., as many as the ceil ($\frac{t}{s}$). This iterative formula allows us to map all non-ascending ordered jumps from the original partitions. If we combine this result with the previous result for the first jumps, we will arrive at the formula.
\begin{equation}
\centering
p(n) = \sum_{i=0}^{floor(\frac{n}{2})} \mathscr{P}(i,n+1-2i)
\end{equation}
\end{proof}

Ono, Folson, Bruinier and Kent  \cite{folsom2012} have shown that there is a fractal-like property in the partition function. The recursive equation that we now show also shows this fractal-like nature in the partitions function. This fractal-like property exists in many other number theory functions as well \cite{hashimoto2007fractal,cattani2016fractal,lapidus1995riemann}. In the next section of this paper, an iterative equation is also derived for the number of divisors function, d(n) and its nature of self-similarity is established using a construction that we will call the network of number of divisors. 

\section{Divisor function d(n)}
Wang Zheng Bing, Robert Fokkink and Wen Fokkink rediscovered an important theorem that relates our earlier work on partitions to the function of the number of divisors \cite{wang1995relation}. Wang and collaborators first establish the definition of a \textit{partition with unequal terms} as a partition that has all its terms unequals, for example, the partition of 7 equal to 4 + 2 + 1 is a partition with unequal terms. A definition that Wang also uses is that of \textit{odd or even partitions}, which are defined as partitions whose \textit{number of terms} is odd or even. Wang defines an \textit{odd or even unequal partition} as a partition that has all its unequal terms and also its number of terms is odd or even, respectively. Based on a different approach previously to the Wang et al. aproach, this theorem was established by D. M. Bressoud and M.V. Subbarao \cite{bressoud1984uchimura}, who base their work on the relationships between partitions and divisors that Uchimura initially proved \cite{uchimura1981identity}.

\begin{thm}[Uchimura-Bressoud-Wang's theorem]
$p_s(n) = d(n)$ for all positive natural numbers n.\\
Where $p_s(n)$ denotes the sum of the \textit{smallest} terms of odd unequal partitions of n minus the \textit{smallest} terms of even unequal partitions of n, and d (n) denotes the number of divisors of n.
\end{thm}

For example, there are five partitions of 7 into unequal parts:\\
\begin{center}
1 + 2 + 4;\quad 1 + 6;\quad 2 + 5;\quad 3 + 4; \quad 7
\end{center}
Since the partitions 1 + 2 + 4 and 7 contain an odd number of summands, they are called odd partitions, whereas the other three partitions are called even. Add the smallest numbers of the odd partitions,1 + 7 = 8, and do the same for the smallest numbers of the even partitions,1 + 2 + 3 = 6. The  difference between these two sums, $p_s(7)$ = 8-6 = 2, is exactly the number of divisors of the prime 7, $d(7)$.\\

Now we will build a special jump sequence different from the previous jump sequences that we have worked on in the section on partitions.\\

\begin{defn}
The transformation of a partition of unequal terms of the type $a_1\; a_2\; \cdots \; a_m$ into a partition of unequal terms such as $a_1\; r\; a_2\; \cdots \; a_m $ is called a \textbf{jump with invariant initial term $a_1$ and order $r$}, with the condition that r is less than any term in the partition and that $a_1$ is not equal to or less than the last term $a_m$. This jump type converts a partition of a number n into a partition of a number n + r.
\end{defn}

\begin{defn}
It is called \textbf{a sequence of descending jumps with an invariant initial term $a_1$ and unequal terms}, the set of jumps with an invariant initial term $a_1$ and whose orders follow a descending sequence.
\end{defn}

Thus, for example, if we start from the number 4. The sequence of descending jumps of invariant initial term 4 and unequal terms that begins with an invariant jump of order 3 and continues with invariant jumps of order 2 and 1 is:\\
4\\
4 3\\
4 2 3\\
4 1 2 3\\

\begin{defn}
\emph{The complete set of descending jumps with an invariant initial term $a_1$ and unequal terms} is the set of all possible sequences of descending jumps with an invariant initial term $a_1$ and unequal terms.
\end{defn}

For example, the complete set of descending jumps of invariant initial term 4 and unequal terms is: \\
4\\
4 1\\
\\
4 2\\
4 1 2\\
\\
4 3\\
4 1 3\\
4 2 3\\
4 1 2 3\\

In another example, the complete set of descending jumps of invariant initial term 5 and unequal terms is: \\
5\\
5 1\\
\\
5 2\\
5 1 2\\
\\
5 3\\ 
5 1 3\\
5 2 3\\
5 1 2 3\\
\\
5 4\\ 
5 1 4\\
5 2 4\\
5 1 2 4\\
5 3 4\\
5 1 3 4\\
5 2 3 4\\
5 1 2 3 4\\

\begin{lem}
All partitions of unequal terms of a number n are contained in a complete set of descending jumps of invariant initial term $a_1$ and of unequal terms, where $a_1$ would be  1, 2 ... or n. 
\end{lem}

\begin{proof}
n is included in the complete set of descending jumps of invariant initial term n. Any other n partition of unequal terms can be ordered in a descending sequence such as $a_{2}\; a_{3}\; \cdots \;  a_{m}\; a_{1}$. If we place the last term as the first term, we have the representation of this partition in the following form $a_{1}\; a_{2}\; a_{3} \; \cdots \;  a_{m}$. The number $a_1$ will be greater than all the other terms in this partition and $a_1$ is less than n. We are going to show that this partition $a_{1}\; a_{2}\; a_{3} \; \cdots \;  a_{m}$ is contained in the set of descending jumps of invariant initial term $a_1$ and unequal terms. Indeed, the terms $a_{2},a_{3},\cdots,a_{m}$ tell us that this partition comes from the successive application of the invariant jumps in $a_1$ of order $a_m, ..., a_3$ and $a_2$. These orders being all different since the partition we are considering has all its different terms. Therefore, this partition will be within the complete set of descending jumps of invariant initial term $a_1$ and unequal terms. 
\end{proof}

\begin{figure}
\centering
\includegraphics[scale=0.3]{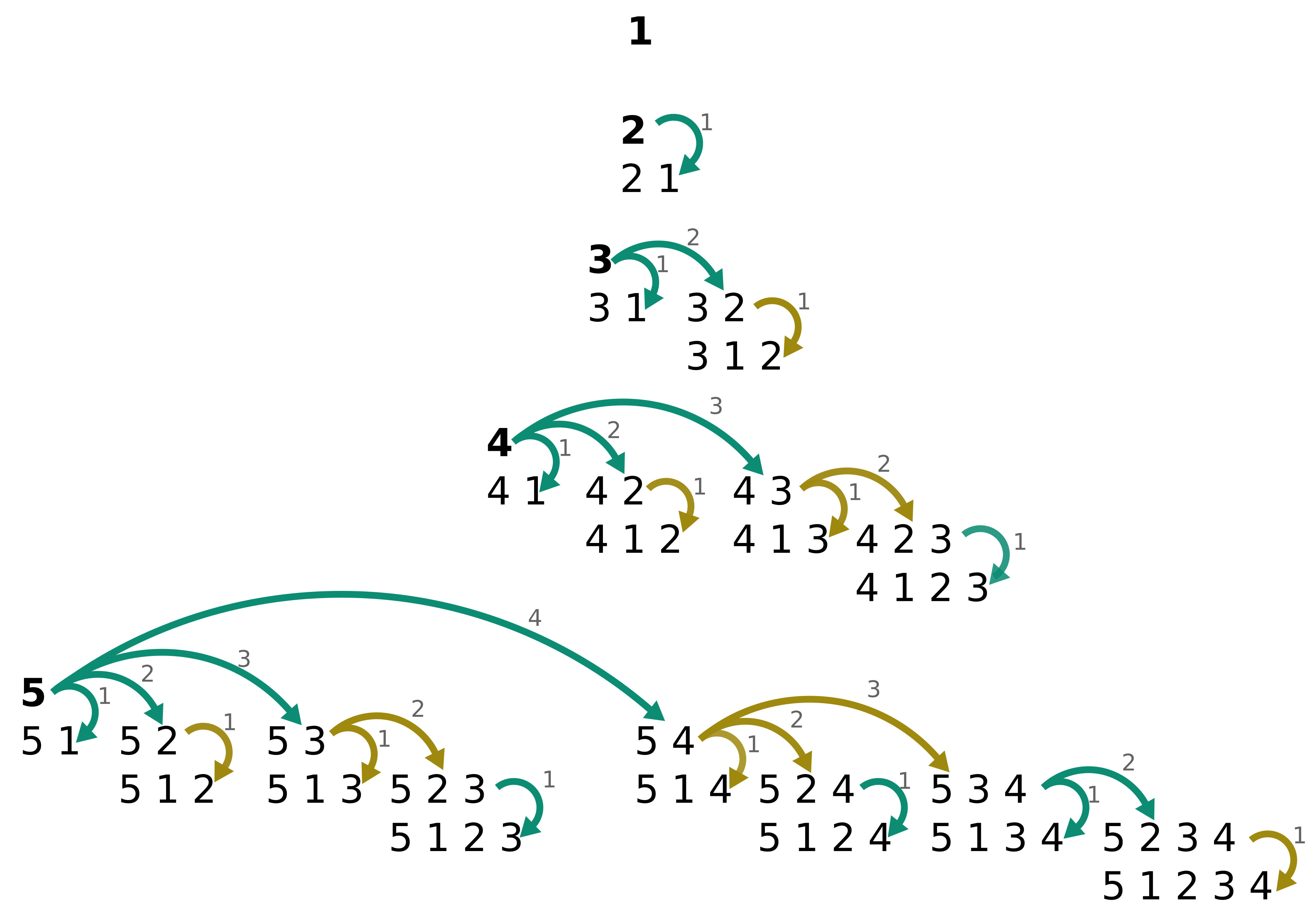}
\caption[]{Complete set of descending jumps of invariant initial term $a_1$ , where $a_1=1,2,3,4$ and $5$ .}
\label{Red_desiguales}
\end{figure}

It is possible to represent a complete set of descending jumps of invariant initial term $a_1$ and unequal terms in a graphical way as shown in Figure \ref{Red_desiguales}. In this representation, the order of the jumps is marked with a gray color. In green, those jumps that are going to generate a partition whose number of terms is even have been indicated, while in brown, those jumps that are going to originate partitions whose number of terms is odd.

\begin{cor}
The partitions of the number n whose terms are all different are those partitions that result from a succession of descending jumps of the invariant initial term $ a_1 $ and of successive orders $ a_m,\cdots, a_3, a_2 $ with the condition that the sum of the jump orders , $ a_m + \cdots + a_3 + a_2$ is equal to  $n - a_1$. 
\end{cor}
\begin{proof}
This result arises from the fact that every partition of the number n is contained within some of the partitions that arise from the sequence of jumps whose invariant initial term $ a_1 $ can take the value from 1 to n. So there is a value of $ a_1 $ such that the successive orders of invariant jumps $ a_m, \cdots, a_3, a_2 $ originates a partition of n, therefore $ a_m + \cdots + a_3 + a_2 + a_1 = n $ y then $ a_m + \cdots + a_3 + a_2 = n - a_1 $
\end{proof}

\begin{lem}
The smallest term of any partition resulting from a succession of descending jumps with invariant initial term is always the order of the last invariant jump performed. 
\end{lem}
\begin{proof}
Since it is a descending succession of invariant jumps, then the order of the successive jumps follows a descending succession, therefore the order of the last jump is the smallest possible term of the resulting partition.
\end{proof}

\begin{figure}
\centering
\includegraphics[scale=0.25]{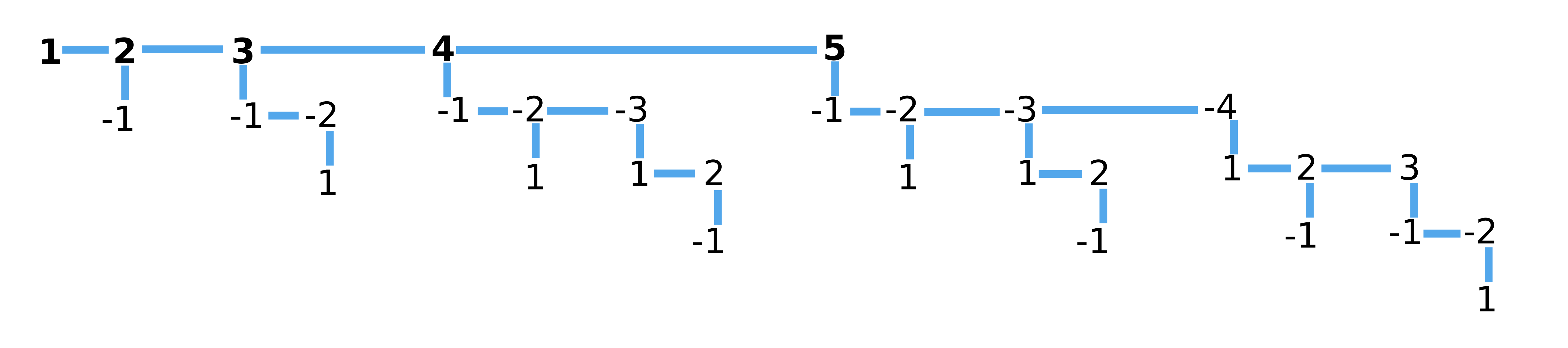}
\caption[]{The smallest terms in the partition that come from a complete set of descending jumps of invariant initial term $a_1$ , where $a_1=1,2,3,4$ and $5$ . If the partition has a even number of terms a minus sign has been added.}
\label{Red_desiguales_smalest}
\end{figure}

In order to relate the representation in the Figure \ref{Red_desiguales} with the Uchimura-Bressoud-Wang theorem, this representation will be changed to the representation of the Figure \ref{Red_desiguales_smalest}, where only the smallest terms of each partition will be placed. Also, in those partitions whose number of terms is even, the minus sign has been added. This representation takes the form of a network that we will call a network of number of divisors, a name that will make sense in the following Lemma.

\begin{defn}
It is called a \textbf{network of number of divisors} to that network whose first vertices are the invariant terms $ a_1 $. The following vertices correspond to the smallest terms of the partitions that result from a succession of descending jumps with an invariant initial term $ a_1 $. Additionally, in those partitions with a number of even terms, the minus sign is added. 
\end{defn}
\begin{defn}
It is called a \textbf{trace of a natural number n}, trace (n), to the sum of the vertices of a network of number of divisors that are the smallest terms of the partitions of n with different terms. 
\end{defn}
\begin{lem}
Trace (n) = d (n), n a natural number and d(n) the number of divisor of n.
\end{lem}
\begin{proof}
This result follows immediately from Uchimura-Bressoud-Wang's Theorem. Indeed, we have previously shown that every partition with different terms is contained in a complete set of descending jumps of invariant initial term $ a_1$ and unequal terms, being able $a_1$ to take the values from 1 to n. Likewise, each partition can only appear once in this complete set because all the terms are different. If we add to this that in a network of number of divisors the partitions have been replaced by their smallest terms and also that the minus sign is added in the partitions that have an even number, then it immediately follows that the trace of the natural number n , as it has been defined, is equal to the number of divisors of the number n.
\end{proof}

\subsection{Theorem in the number of divisor of a natural number}
\begin{thm}
If d(n) is the divisor function and if we define a vector $\overline{E_n}$ as :\\
\begin{equation}
\begin{split}
\overline{E_n}\;  = &1\,A_{\tiny {n\,n-1}} A_{\tiny{n-1\,n-2}} \cdots A_{\tiny{3\,2}} A_{\tiny{2\,1}} \overline{I_1} +\\
                    &2\,A_{\tiny {n\,n-1}} A_{\tiny{n-1\,n-2}} \cdots A_{\tiny{3\,2}}\overline{I_2} + \\
                    &\cdots +\\
                    &(n-1)\,A_{\tiny {n\,¨n-1}} \overline{I_{n-1}} +\\
                    &n\,\overline{I_n}
\end{split}
\end{equation}
Then:
\begin{equation}
\overline{E_n}[i] = d(i),\;i = t(n-1)+1,\cdots ,t(n)
\end{equation}
Where,
$\overline{I_i} [k]=
  \begin{cases}
   1, & if \; k = t(i-1)+1 \\
   0, & if \; k \ne t(i-1)+1
  \end{cases}$\\
$
A_{\tiny{i\,j}} [k,l]= 
  \begin{cases}
   -1, & if \; k = l \\
	1, & if \; k-i = l  \\
   0, & if \; other \; cases
  \end{cases}\\
$
With $k = 1,\cdots,t(i)$; $l = 1,\cdots,t(j)$. t(n) =  nth tringular number. Also, the size of a matrix $A_{\tiny{i\,j}} $ is t(i) x t(j).
\end{thm}
\begin{proof}
In order to proof this theorem, the Figure \ref{Red_divisores_grupos} can help us in order to understand the generation of a special vector. First, we note that the network of divisors is a fractal network in the sense that a part of it is repeated for the network to evolve. We have placed the repeating parts in circles of the same color. 
\begin{figure}[H]
\centering
\includegraphics[scale=0.25]{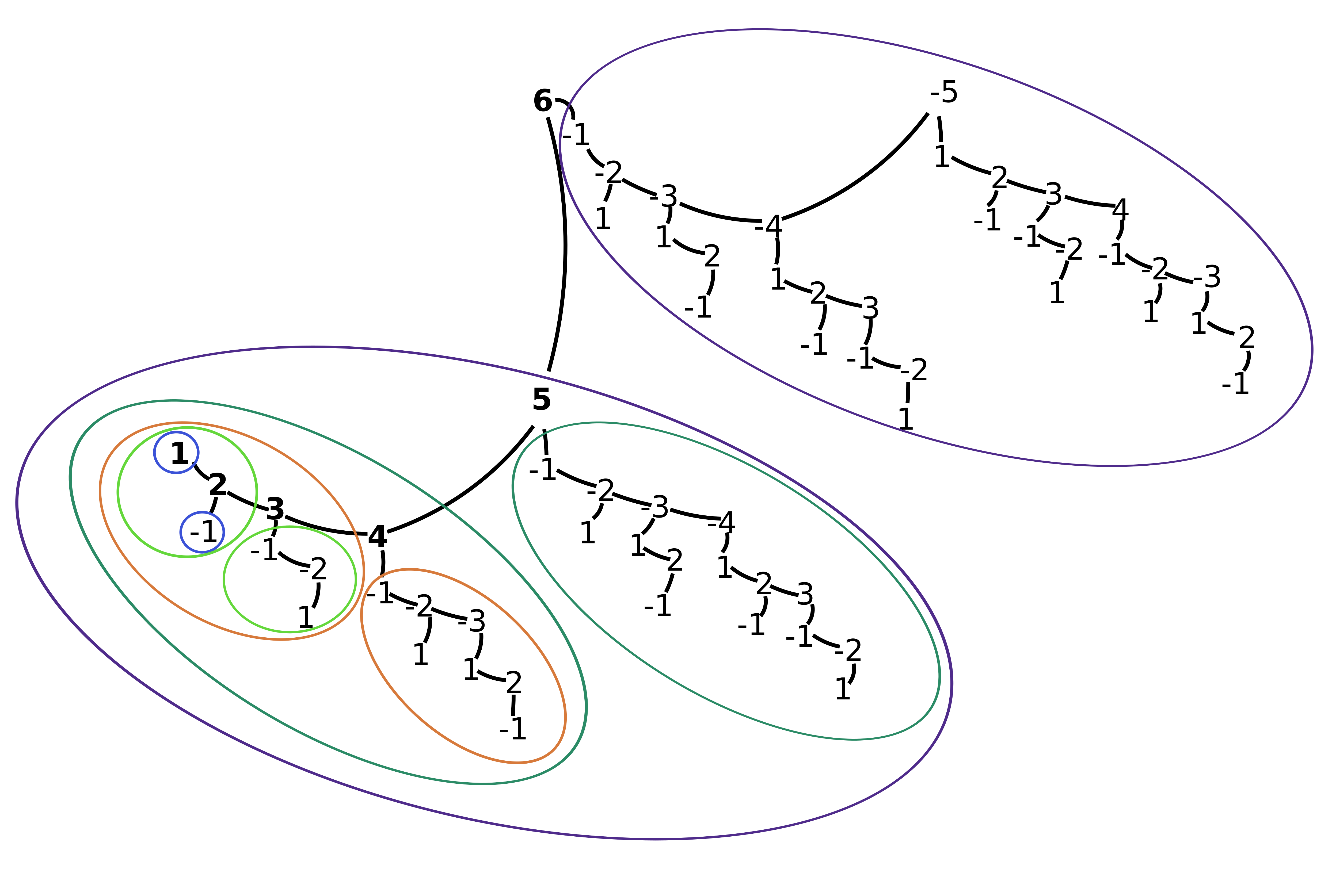}
\caption[]{The groups in the divisor network until number 6 shows that this network is really a fractal network.}
\label{Red_divisores_grupos}
\end{figure}
The only difference between these parts with the same color is the sign of the elements that make it up. On the other hand, the trace of a given number n is the sum of all those numbers that appear in the network of divisors that come from the partitions of n. It is very easy to locate them since reaching these numbers from number 1, located at the beginning of the network, takes n steps (understanding steps such as visiting a certain number on the network following the network's plot). For example:\\ 
trace (6) = 6 + (-1) + (-2) + 1 = 4 = d (6)\\
trace (5) = 5 + (-1) + (-2) = 2 = d (5)\\
trace (4) = 4 + (-1) = 3 = d (4)\\
trace (3) = 3 + (-1) = 2 = d(3)\\
trace (2) = 2 = d(2)\\
trace (1) = 1 = d(1) \\
\begin{figure}[H]
\centering
\includegraphics[scale=0.14]{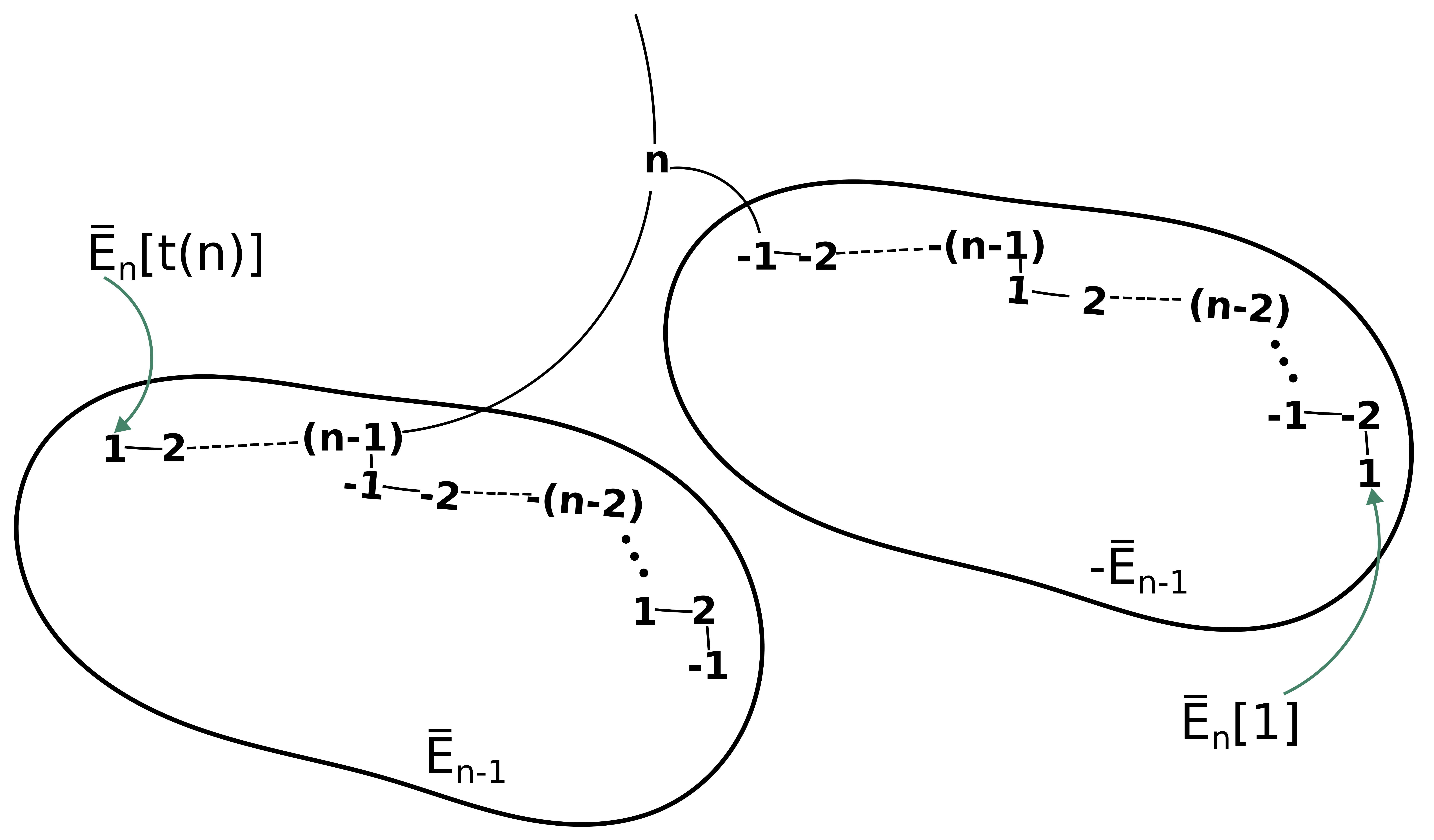}
\caption[]{Abstract representation of a divisor network.}
\label{Red_divisores_abstracto}
\end{figure}
From the Figure \ref{Red_divisores_grupos}, we are going to derive a more general Figure, the Figure \ref{Red_divisores_abstracto}, which relates the part of the network before the number n, and the part of the network that originates from the number n. As seen in the Figure, these two parts are identical only differing in the sign of their elements. We will call the first part the initial block and the second part the final block. We will define the vector $\overline{E_n}$ in this general network. This vector will be generated by adding the terms of the network in the same direction in which the trace of any number in the network is calculated. To generate this vector, only the initial and final blocks related to the number n are used. For example:\\
$
\overline{E_1} =
\begin{bmatrix}
\textbf{1}
\end{bmatrix}
$;\;
$
\overline{E_2} =
\begin{bmatrix}
-1 \\
\textbf{2} \\
1 
\end{bmatrix} 
$;\;
$
\overline{E_3} =
\begin{bmatrix}
1  \\
-2 \\
-1 \\
\textbf{2 } \\
2  \\
1 
\end{bmatrix} 
$;\;
$
\overline{E_4} =
\begin{bmatrix}
-1 \\
2  \\
1  \\
-2 \\
-1  \\
-3 \\
\textbf{3} \\
2  \\
2  \\
1 
\end{bmatrix} 
$;\;
$
\overline{E_5} =
\begin{bmatrix}
1  \\
-2 \\
-1 \\
2  \\
1  \\
2 \\
-1  \\
-1  \\
-4 \\
-2 \\
\textbf{2} \\
3 \\
2  \\
2  \\
1 
\end{bmatrix} 
$\\
In general, each vector $\overline{E_n}$ has a size of t(n), where t(n) is the nth triangular number. We can extract from the Figure \ref{Red_divisores_abstracto} that:\\
$\overline{E_n} = \begin{bmatrix}
0\\
\vdots\\
0\\
\vdots\\
\overline{E_{n-1}}
\end{bmatrix} $ + $\begin{bmatrix}
0\\
\vdots\\
n\\
\vdots\\
0
\end{bmatrix} 
$ + $\begin{bmatrix}
-\overline{E_{n-1}}\\
\vdots\\
0\\
\vdots\\
0
\end{bmatrix} 
$ = $\begin{bmatrix}
-\overline{E_{n-1}}\\
\vdots\\
0\\
\vdots\\
\overline{E_{n-1}}
\end{bmatrix} $ + $\begin{bmatrix}
0\\
\vdots\\
n\\
\vdots\\
0
\end{bmatrix} 
$ \\
where the number n in the vector $\begin{bmatrix}
0\\
\vdots\\
n\\
\vdots\\
0
\end{bmatrix} 
$  it occupies the (t(n-1)+1)-th position. Also, $\overline{E_n}[i] = d(i)$ , only for $t(n-1)+1,\cdots ,t(n)$. From this equation we deduce that:\\
\begin{equation}
\overline{E_n} = A_{\tiny {n\,n-1}} \overline{E_{n-1}} + n\overline{I_n}
\end{equation}
Applying this equation iteratively to $\overline{E_{n-1}} , \cdots, \overline{E_2}$ and $\overline{E_1}$ , we obtain:\\
$\overline{E_n}\;  = \; 1\,A_{\tiny {n\,n-1}} A_{\tiny{n-1\,n-2}} \cdots A_{\tiny{3\,2}} A_{\tiny{2\,1}} \overline{I_1} + 2\,A_{\tiny {n\,n-1}} A_{\tiny{n-1\,n-2}} \cdots A_{\tiny{3\,2}}\overline{I_2} + \cdots + (n-1)\,A_{\tiny {n\,¨n-1}} \overline{I_{n-1}} + n\,\overline{I_n}$ \\
Where it is shown that $\overline{E_n}[i] = d(i)$ , only for $i = t(n-1)+1,\cdots ,t(n)$.\\
Where,
$\overline{I_i} [k]=
  \begin{cases}
   1, & if \; k = t(i-1)+1 \\
   0, & if \; k \ne t(i-1)+1
  \end{cases}$\\
$
A_{\tiny{i\,j}} [k,l]= 
  \begin{cases}
   -1, & if \; k = l \\
	1, & if \; k-i = l  \\
   0, & if \; other \; cases
  \end{cases}\\
$
With $k = 1,\cdots,t(i)$; $l = 1,\cdots,t(j)$. t(n) =  nth tringular number

\end{proof}

\section{Concluding Remarks}
Thanks to the kind suggestion of Robbert Fokkink, who in turn received this suggestion from George E. Andrews, I had the opportunity to learn about the work of the Bressoud and Subbarao. Bressoud and Subbarao establish a general relationship between the functions $ \sigma_i (n) $ and the partitions of unequal terms, where $\sigma_i (n)$ is the sum of the ith powers of the divisors of n. This equation for the case of the function $ \sigma_1(n)$, sum of divisors of the number n, is established by:
\begin{thm}[Bressoud-Subbarao's theorem]
\begin{equation}
\sigma_1(n)  = -\sideset{}{'}\sum_{\tiny \pi \vdash n}^{} (-1)^{\#\pi} \sum_{\tiny j = 1}^{\tiny \lambda(\pi)} (L(\pi)-\lambda(\pi) +j)
\end{equation}
\end{thm}
Where, $\pi \vdash n $ means that $\pi$ is a partition of n, the prime ($'$) on the summation restricts the sum to those partitions which have
distinct parts, $\#(\pi)$ is the number of terms in $\pi$ and $\lambda(\pi)$ is the smallest term in $\pi$ and $L(\pi)$ is the largest term in $\pi$. Using this equation in Figure \ref{Red_desiguales} we can construct a network that we will call a network of sum of divisors, analogous to the network of number of divisors. The trace, defined in a general way as those terms of the network that are at the same distance n from the initial point 1, gives us the sum of the divisors of the number n. 
\begin{figure}[H]
\centering
\includegraphics[scale=0.25]{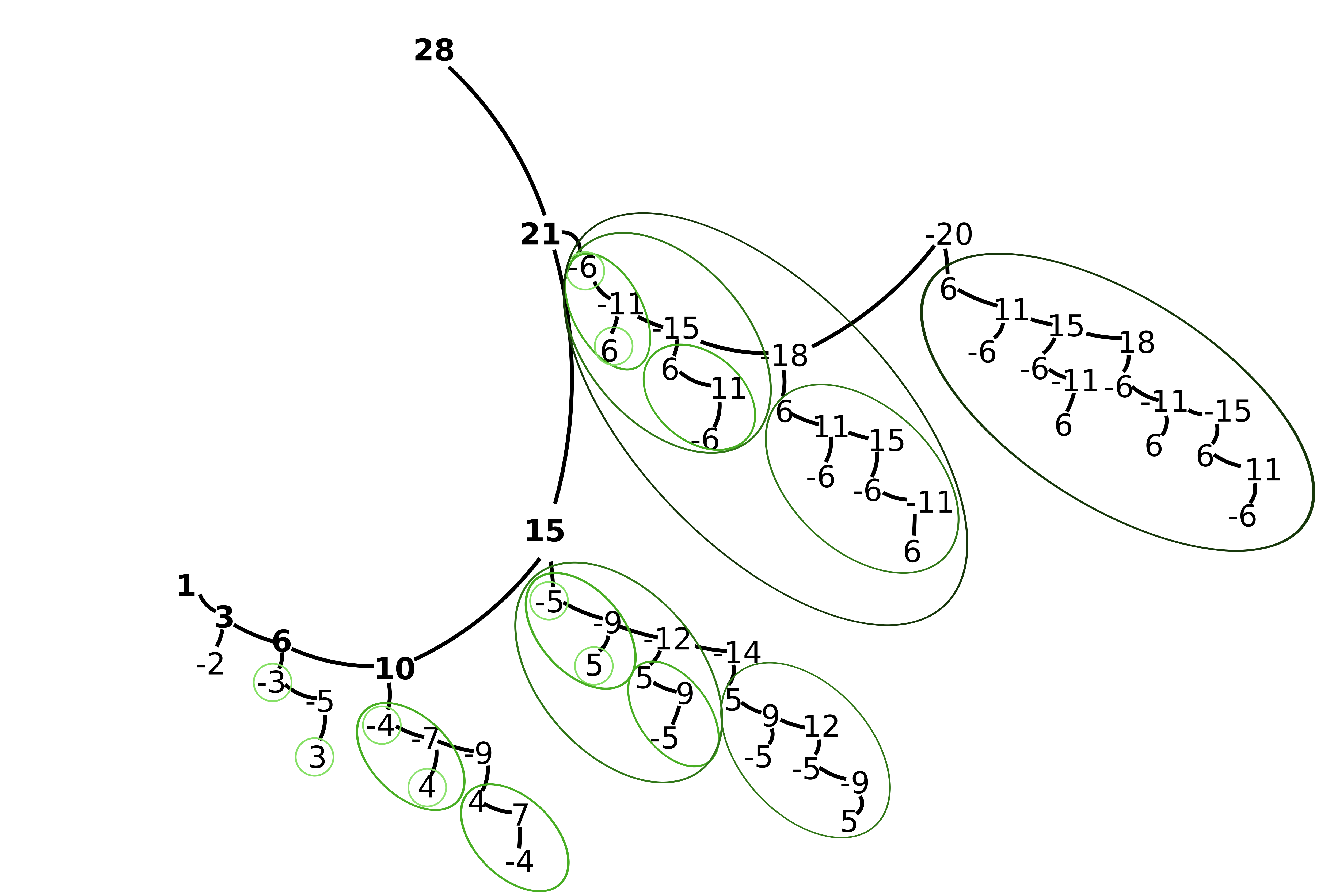}
\caption[]{The groups in the sum of divisor network until number 7 shows that this network also a special self-similarity network.}
\label{Red_suma-divisores}
\end{figure}
Using this general definition of trace mentioned above, and calculating it in Figure for a number between 1 and 7, we can verify that:\\
trace (1) = 1 = $\sigma_1(1)$,\\
trace(2) = 3 = $\sigma_1(2)$,\\
trace(3) = 6 + (-2) = 4 = $\sigma_1(3)$,\\
trace(4) = 10 + (-3) = 7 = $\sigma_1(4)$,\\
trace(5) = 15 + (-4) + (-5) = 6 = $\sigma_1(5)$,\\
trace(6) = 21 + (-5) + (-7) + 3 = 12 = $\sigma_1(6)$,\\
trace(7) = 28 + (-6) + (-9) + (-9) + 4 = 8 = $\sigma_1(7)$,\\

This network of sum of divisors shows a main branch composed of the triangular numbers of the number to which you want to calculate the sum of divisors. From each of these numbers a fractal-type network arises that begins with the number -n and evolves in a slightly different way than the number of divisors network evolves, but maintains its self-similarity behavior. The present work thus shows that the networks of divisors d(n) and $ \sigma(n) $ have a structure of self similarity.

\section*{Acknowledgement}
The author want to express their sincerest gratitude to FONDECYT (Convenio nº208-2015-FONDECYT) for his Master scholarship. Likewise, express its gratitude to Dr. Jim Cleaves from the Blue Marble Space Institute of Science (BMSIS) with whom he had the opportunity to learn about the networks of chemical reactions where he read about integer hyperflows and where the idea of looking for an equation for integer partitions came up in order to use that equation in future chemistry works. The author would like to especially thank Dr. Robert Fokkink for showing him the work of Bressoud-Subbarao, where the equation on which the generation of the network of sum of divisors is based appears. The author would like to thank his thesis advisor in Labinfis, Dr. Gino Picasso, for his encouragement to allow the author to move away to other topic of interes and future importance. Likewise, he thanked his love Ruth Quispe, his parents Julia Simbron and Antonio Cruz, and family for their incessant support without whose intervention this work would not have been done. Finally, he would like to thank his past teachers (university and collage), friends in general who have been reviewing this draft and giving me their criticism and encouragement. 

\bibliographystyle{apacite}
\bibliography{referencias/References}

\end{document}